\documentclass[12pt]{myarticle}

\newtheorem{thm}{Theorem}

\newtheorem{cor}{Corollary}

\usepackage{tikz}
\usetikzlibrary{positioning}
\usepackage{subfigure}
\usepackage{caption}
\usepackage{graphicx}
\usepackage{pgfkeys}
\usepackage{bm}
\usepackage[all]{xy}

\usepackage{amsmath}

\usepackage{txfonts}
\usepackage[colorlinks]{hyperref}

\begin{document}

\title{The general Zagreb index of lattice networks}

\author[ps]{Prosanta SARKAR}
\ead{prosantasarkar87@gmail.com}

\address[ps]{Department of Mathematics, National Institute of Technology, Durgapur, India.}

\author[nd]{Nilanjan DE}
\ead{de.nilanjan@rediffmail.com}

\address[nd]{Department of Basic Sciences and Humanities (Mathematics),\\ Calcutta Institute of Engineering and Management, Kolkata, India.}


\author[ap]{Anita PAL}
\ead{anita.buie@gmail.com}

\address[ap]{Department of Mathematics, National Institute of Technology, Durgapur, India.}

\begin{abstract}
A topological index is a real number which is derived from a network or a graph by mathematically that characterizes the whole of its structural properties. Recently, there are various topological indices that have been introduced in mathematical chemistry to predict the properties of molecular topology. Among, the degree based topological indices such as Zagreb indices, forgotten topological index, redefined Zagreb index, Randi\'{c} index, general first Zagreb index, symmetric division deg index and hence so forth are most important, because of their chemical significance. In this work, we study the general Zagreb index of hexagonal and triangular lattice networks.
  
\medskip

\end{abstract}


\maketitle

\section{Introduction}

Networks are mathematical structures used to study pairwise relationships between objects and entities. It has found many applications in chemical, physical and social sciences, information technology, and optimization. In network research, it relies on that finding a
suitable measure and use this measure to quantify network robustness. A network or a graph is a set of nodes or vertices connected by a set of lines or edges. More generally, a network or a graph is an ordered pair of two sets namely vertex set $V(G)$ and the edge set $E(G)$ that is $G=(V(G), E(G))$. The degree of a vertex $v\in V(G)$ is the number of adjacent vertices to $v$ and is denoted by $d_G(v)$. A graph invariant or a topological index is a numeric quantity associated with a graph that characterizes the whole of its topology. In chemical sciences, there are various types of topological indices have been introduced such as degree based topological indices, distance-based topological indices, spectrum-based topological indices and so forth. Among the vertex degree based topological indices such as Zagreb indices, forgotten topological indices, redefined Zagreb index, general first Zagreb index, general Randi\'{c} index and symmetric division deg index are most studied and widely used in chemical sciences to predict structural properties of molecules. To compute the total $\pi$-electron energy$(\epsilon)$ of carbon atoms Gutman and Trinajesti\'{c} \cite{gutm72}, were first introduced the Zagreb indices in 1972 and are defined as
  \[{{M}_{1}}(G)=\sum\limits_{v\in V(G)}{{{d}_{G}}{{(v)}^{2}}}=\sum\limits_{uv\in E(G)}{[{{d}_{G}}(u)+{{d}_{G}}(v)]}\] and \[{{M}_{2}}(G)=\sum\limits_{uv\in E(G)}{{{d}_{G}}(u){{d}_{G}}(v)}.\] We refer our reader to \cite{gut18,ndep17,sidd16}, for further study about this index. The ``forgotten topological index'' was discovered in the same paper where Zagreb indices were studied first time in the paper \cite{gutm72}, in 1972 but then not much studied about this index, Furtula and Gutman reinvestigate this index again in a paper \cite{furt15}, in 2015 and is defined as \[{F}(G)=\sum\limits_{v\in V(G)}{{{d}_{G}}{{(v)}^{3}}}=\sum\limits_{uv\in E(G)}{[{{d}_{G}}(u)^2+{{d}_{G}}(v)^2]}.\] See \cite{den17,ram16,akh17,thay18}, for some recent study about this index. In \cite{ranj13}, Ranjini et al. was first introduced the redefined Zagreb index in 2013, and is defined as  
\[{{ReZM}}(G)=\sum\limits_{uv\in E(G)}{{{d}_{G}}(u){{d}_{G}}(v)[{{d}_{G}}(u)+{{d}_{G}}(v)]}.\] We refer our reader to \cite{bzha16,gao16,kum17}, for some recent study about this index. Based on first Zagreb index and F-index Li and Zheng in \cite{li05} was first introduced the general Zagreb index and is defined as follows \[M^{\alpha}(G)=\sum\limits_{u\in V(G)}d_G(u)^{\alpha}\] where, $\alpha \neq 0,~1~and ~\alpha \in {\rm I\!R}$. Interested reader to see \cite{liu10,bed18,xav14}, for further study about this index. Gutman and Lepovi\'{c} generalized the Randi\'{c} index in \cite{gut01}, and is defined as
 \[R_{a}=\sum\limits_{uv\in E(G)}\{d_G(u).d_G(v)\}^{a}\] where, $a \neq 0,~a \in {\rm I\!R}$. For further study about this index, we refer our reader to \cite{zho09,cui17}. The symmetric division deg index is defined as \[SDD(G)=\sum\limits_{uv\in E(G)}[\frac{d_G(u)}{d_G(v)}+\frac{d_G(v)}{d_G(u)}].\] See \cite{alex14,gup17}, for further study about this index. Flowed by Zagreb indices Azari et al.\cite{azar11}, were first introduced the general Zagreb index or $(a,b)$-Zagreb index in 2011 and is defined as \[Z_{a,b}(G)=\sum\limits_{uv\in E(G)}(d_G(u)^ad_G(v)^b+d_G(u)^bd_G(v)^a).\] For some recent study about this index we refer our reader to \cite{sar18,sar19,sark18}.  
\begin{table}[ht]
\caption{Relations between (a,b)-Zagreb index with some other topological indices:} 
\centering
\begin{tabular}{|c|c|}
\hline
Topological index & Corresponding (a,b)-Zagreb index \\
\hline
First Zagreb index ${{M}_{1}}(G)$ & $Z_{1,0}(G)$ \\
\hline
Second Zagreb index ${{M}_{2}}(G)$ & $\frac{1}{2}Z_{1,1}(G)$ \\
\hline
Forgotten topological index ${F}(G)$ & $Z_{2,0}(G)$\\
\hline
Redefined Zagreb index ${{ReZM}}(G)$ & $Z_{2,1}(G)$\\
\hline
General first Zagreb index $M^a(G)$ & $Z_{a-1,0}(G)$\\
\hline
General Randi\'{c} index $R_{a}$ & $\frac{1}{2} Z_{a,a}$\\
\hline
Symmetric division deg index $SDD(G)$ & $Z_{1,-1}(G)$\\
\hline
\end{tabular}
\label{table:nonlin} 
\end{table}
In this paper, we study the general Zagreb index of lattice networks. A lattice network or simply a lattice is set of all linear combinations of vectors in ${\rm I\!R^n}.$ Lattices have some interesting features that make them important agents for used in matter physics. For the symmetric nature of its topology, lattice graphs are the most common network structures and are widely used in distributed parallel computation, distributed control, wired circuits such as CMOS and CCD based devices, satellite constellations and sensor networks. In addition, the improvement of micro and nano-technologies has also empowered the deployment of sensor networks for estimating and observing purposes. For the same reason, the study of various physical and chemical indices of various lattice graphs has captured the great attention to many researchers. Let us suppose that $G_1$ and $G_2$ be two graphs with vertex sets $V(G_1),~V(G_2)$ and edge sets $E(G_1),~E(G_2)$. The Cartesian product between $G_1$ and $G_2$ is defined as $G_1 \square G_2$, and any two vertices $(a,a^{\prime})$ and $(b,b^{\prime})$ are connected in $G_1 \square G_2$ if and only if either $a=b$ and $a^{\prime}b^{\prime} \in E(G_2)$ or $a^{\prime}=b^{\prime}$ and $ab\in E(G_1)$ where $a,b\in V(G_1)$ and $a^{\prime},b^{\prime}\in V(G_2)$. A triangular lattice network of $m\times n$ with toroidal boundary condition is nothing but the Cartesian product between $C_m$ and $C_n$ $(m,~n \ge 3)$ with an additional diagonal edge added to every square, as shown in Figure 6 and is defined as $T^t[m,n]$. The triangular lattice network with cylindrical boundary condition $T^c[m,n]$ is derived from $T^t[m,n]$ by removing the edges $(a_1,a_1^*),~(a_2,a_2^*),.......,(a_n,a_n^*)$; $(a_2,a_1^*),~(a_3,a_2^*),.......,(a_n,a_{n-1}^*),~(a_1,a_n^*)$, as shown in Figure 5. The triangular lattice with free boundary condition $T^f[m,n]$ is obtained from $T^c[m,n]$ by further delete the edges $(b_1,b_1^*),~(b_2,b_2^*),.......,(b_m,b_m^*)$; $(b_2,b_1^*),~(b_3,b_2^*),.......,(b_m,b_{m-1}^*)$, as shown in Figure 4. The hexagonal lattice is lot like a Cartesian lattice (square) except pixels in hex lattice are hexagons instead of squares as it is the case for square lattice. The hexagonal lattice with toroidal boundary condition is defined as $H^t[m,n]$, and shown in Figure 3. The hexagonal lattice with cylindrical boundary condition is derived from $H^t[m,n]$ by removing the edges $(b_1,a_1^*),~(a_1,a_2^*),......,(a_{n-1},a_n^*),~(a_n,c_{m+1})$, as shown in Figure 2 and is defined as $H^c[m,n]$. If further we delete the edges $(b_1,c_1),~(b_2,c_2),.......,(b_{m+1},c_{m+1})$ from $H^c[m,n]$, then we obtained hexagonal lattice network with free boundary condition as shown in Figure 1 and is defined as $H^f[m,n]$. Followed by lattice networks, J.B. Liu and X.F. Pan studied the asymptotic incidence energy of lattices in \cite{liu15} and L. Chen et al. studied the Wiener polarity index of lattice networks in \cite{chen16}.  In this paper, we compute the general Zagreb index of hexagonal and triangular lattice networks with toroidal, cylindrical and free boundary conditions.  

\section{Main Results:}

In the following section we compute the general Zagreb index or the $(a,b)$-Zagreb index of hexagonal and triangular lattice networks with toroidal, cylindrical and free boundary conditions. At first, we start with hexagonal lattice network with free boundary condition $H^f[m,n].$ Note that the total number of vertices of $H^f[m,n]$ or $|V(H^f[m,n])|=2(m+1)(n+1)$ and $|E(H^f[m,n])|= m(n+1)+(m+1)(2n+1)$. Based on degree of the end vertices of each edge of $H^f[m,n]$ we divide the edge set of $H^f[m,n]$ into four sub sets namely $E_1(H^f[m,n])$, $E_2(H^f[m,n])$, $E_3(H^f[m,n])$, and $E_4(H^f[m,n])$ as shown in the table 2.
\begin{table}[ht]
\caption{Edge partition of hexagonal lattice network with free boundary condition} 
\centering  
\begin{tabular}{|c|c|}
\hline
$(d(u),d(v)): ~uv\in E(H^f[m,n])$ & Total number of edges\\
\hline
(3,2) & $4m+4n-4$\\
\hline
(3,3) & $m(n-1)+(2n-1)(m-1)$\\
\hline
(1,3) & $2$\\
\hline
(2,2) & $2$\\
\hline
\end{tabular}
\label{table:nonlin} 
\end{table}
\begin{figure*}[ht]
\centering   
\captionsetup{singlelinecheck = true}
\begin{tikzpicture}[scale=.8,
every edge/.style = {draw=black, thick},
 vrtx/.style args = {#1/#2}{
      circle, draw, thick, fill=black,
      scale=.1, label=#1:#2}
                    ]
\node (A) [vrtx=/]   [label=left:{$b_1$}]   at (0,0) {};
\node (B) [vrtx=/]   [label=below left:{$a_1$}]  at (2,0) {};
\node (C) [vrtx=/]    [label=below left:{$a_2$}] at ( 4, 0) {};
\node (D) [vrtx=/]      at (6, 0) {};
\node (E) [vrtx=/]   [label=below left:{$a_{n}$}]  at ( 8, 0) {};
\node (F) [vrtx=/]    at ( 1, 0.5) {};
\node (G) [vrtx=/]     at ( 3, 0.5) {};
\node (H) [vrtx=/]     at ( 5, 0.5) {};
\node (I) [vrtx=/]    at ( 7, 0.5) {};
\node (J) [vrtx=/]   [label=below right:{$c_1$}] at ( 9, 0.5) {};
\node (K) [vrtx=/]    [label=left:{$b_2$}] at ( 1, 1.5) {};
\node (L) [vrtx=/]    at ( 3, 1.5) {};
\node (M) [vrtx=/]    at ( 5, 1.5) {};
\node (N) [vrtx=/]     at ( 7, 1.5) {};
\node (O) [vrtx=/]     at ( 9, 1.5) {};
\node (P) [vrtx=/]     at ( 2, 2) {};
\node (Q) [vrtx=/]     at ( 4, 2) {};
\node (R) [vrtx=/]     at ( 6, 2) {};
\node (S) [vrtx=/]      at ( 8, 2) {};
\node (T) [vrtx=/]    [label=below right:{$c_2$}] at (10,2) {};
\node (U) [vrtx=/]    [label=left:{$b_{3}$}] at (2,3) {};
\node (V) [vrtx=/]     at ( 4, 3) {};
\node (W) [vrtx=/]     at (6, 3) {};
\node (X) [vrtx=/]     at ( 8, 3) {};
\node (Y) [vrtx=/]     at ( 10, 3) {};
\node (Z) [vrtx=/]     at ( 3, 3.5) {};
\node (A*) [vrtx=/]     at ( 5, 3.5) {};
\node (B*) [vrtx=/]    at ( 7, 3.5) {};
\node (C*) [vrtx=/]    at ( 9, 3.5) {};
\node (D*) [vrtx=/]    at ( 11, 3.5) {};
\node (E*) [vrtx=/]     at ( 3, 4.5) {};
\node (F*) [vrtx=/]      at ( 5, 4.5) {};
\node (G*) [vrtx=/]     at ( 7, 4.5) {};
\node (H*) [vrtx=/]     at ( 9, 4.5) {};
\node (I*) [vrtx=/]     at ( 11, 4.5) {};
\node (J*) [vrtx=/]     at ( 4, 5) {};
\node (K*) [vrtx=/]       at ( 6, 5) {};
\node (L*) [vrtx=/]     at ( 8, 5) {};
\node (M*) [vrtx=/]     at ( 10, 5) {};
\node (N*) [vrtx=/]     [label=below right:{$c_{m}$}] at ( 12, 5) {};
\node (O*) [vrtx=/]     [label= left:{$b_{m+1}$}] at ( 4, 6) {};
\node (P*) [vrtx=/]     at ( 6, 6) {};
\node (Q*) [vrtx=/]     at ( 8, 6) {};
\node (R*) [vrtx=/]     at ( 10, 6) {};
\node (S*) [vrtx=/]     at ( 12, 6) {};
\node (T*) [vrtx=/]    [label=above left:{$a_{1}^*$}] at ( 5, 6.5) {};
\node (U*) [vrtx=/]     [label=above left:{$a_{2}^*$}] at ( 7, 6.5) {};
\node (V*) [vrtx=/]     at ( 9, 6.5) {};
\node (W*) [vrtx=/]     [label=above left:{$a_{n}^*$}] at ( 11, 6.5) {};
\node (X*) [vrtx=/]     [label=below right:{$c_{m+1}$}] at ( 13, 6.5) {};
\node (Y*) [vrtx=/]   at (5.9,1) {};
\node (Z*) [vrtx=/]   at (6,1) {};
\node (A**) [vrtx=/]   at (6.1,1) {};
\node (B**) [vrtx=/]   at (6.9,2.5) {};
\node (C**) [vrtx=/]   at (7,2.5) {};
\node (D**) [vrtx=/]   at (7.1,2.5) {};
\node (E**) [vrtx=/]   at (8.9,5.5) {};
\node (F**) [vrtx=/]   at (9,5.5) {};
\node (G**) [vrtx=/]   at (9.1,5.5) {};
\node (C**) [vrtx=/]   at (3.9,3.9) {};
\node (D**) [vrtx=/]   at (4,4) {};
\node (E**) [vrtx=/]   at (4.1,4.1) {};
\node (C**) [vrtx=/]   at (5.9,3.9) {};
\node (D**) [vrtx=/]   at (6,4) {};
\node (E**) [vrtx=/]   at (6.1,4.1) {};
\node (C**) [vrtx=/]   at (7.9,3.9) {};
\node (D**) [vrtx=/]   at (8,4) {};
\node (E**) [vrtx=/]   at (8.1,4.1) {};
\node (C**) [vrtx=/]   at (9.9,3.9) {};
\node (D**) [vrtx=/]   at (10,4) {};
\node (E**) [vrtx=/]   at (10.1,4.1) {};
\path   (A) edge (F)
        (B) edge (F)
        (B) edge (G)
        (G) edge (C)
        (C) edge (H)
        (H) edge (D)
        (D) edge (I)
        (I) edge (E)
        (E) edge (J)
        (F) edge (K)
        (G) edge (L)
        (H) edge (M)
        (I) edge (N)
        (J) edge (O)
        (K) edge (P)
        (P) edge (L)
        (L) edge (Q)
        (Q) edge (M)
        (M) edge (R)
        (R) edge (N)
        (N) edge (S)
        (S) edge (O)
        (O) edge (T)
        (P) edge (U)
        (Q) edge (V)
        (R) edge (W)
        (S) edge (X)
        (T) edge (Y)
        (U) edge (Z)
        (Z) edge (V)
        (V) edge (A*)
        (A*) edge (W)
        (W) edge (B*)
        (B*) edge (X)
        (X) edge (C*)
        (C*) edge (Y)
        (Y) edge (D*)
        (Z) edge (E*)
        (A*) edge (F*)
        (B*) edge (G*)
        (C*) edge (H*)
        (D*) edge (I*)
        (E*) edge (J*)
        (J*) edge (F*)
        (F*) edge (K*)
        (K*) edge (G*)
        (G*) edge (L*)
        (L*) edge (H*)
        (H*) edge (M*)
        (M*) edge (I*)
        (I*) edge (N*)
        (J*) edge (O*)
        (K*) edge (P*)
        (L*) edge (Q*)
        (M*) edge (R*)
        (N*) edge (S*)
        (O*) edge (T*)
        (T*) edge (P*)
        (P*) edge (U*)
        (U*) edge (Q*)
        (Q*) edge (V*)
        (V*) edge (R*)
        (R*) edge (W*)
        (W*) edge (S*)
        (S*) edge (X*);
\end{tikzpicture}
\caption*{~$H^f[m,n]$}
\captionsetup{singlelinecheck = true, justification=justified}
\caption{The example of hexagonal lattice with free boundary conditions.}
\end{figure*}

\begin{thm}
For $H^f[m,n]$ the general Zagreb index or $(a,b)$-Zagreb index is given by
\begin{eqnarray}
\nonumber Z_{a,b}(H^f[m,n])&=&(3mn-2m-2n+1)(2.3^{(a+b)})+2^{(a+b+2)}\\
&&+(4m+4n-4)(3^a.2^b+3^b.2^a)+2(3^a+3^b).
\end{eqnarray}  
\end{thm}
\textit{Proof.} From definition of general Zagreb index we get, 
\begin{eqnarray*}
Z_{a,b}(H^f[m,n])&=&\sum\limits_{uv\in E(H^f[m,n])}(d_{H^f[m,n]}(u)^ad_{H^f[m,n]}(v)^b+d_{H^f[m,n]}(u)^bd_{H^f[m,n]}(v)^a)\\
&=&\sum\limits_{uv\in E_1(H^f[m,n])}(3^a.2^b+3^b.2^a)+\sum\limits_{uv\in E_2(H^f[m,n])}(3^a.3^b+3^b.3^a)\\
&&+\sum\limits_{uv\in E_3(H^f[m,n])}(1^a.3^b+1^b.3^a)+\sum\limits_{uv\in E_4(H^f[m,n])}(2^a.2^b+2^b.2^a)\\
&=&|E_1(H^f[m,n])|(3^a.2^b+3^b.2^a)+|E_2(H^f[m,n])|(3^a.3^b+3^b.3^a)\\
&&+|E_3(H^f[m,n])|(1^a.3^b+1^b.3^a)+|E_4(H^f[m,n])|(2^a.2^b+2^b.2^a)\\
&=&2(1^a.3^b+1^b.3^a)+[m(n-1)+(2n-1)(m-1)](3^a.3^b+3^b.3^a)\\
&&+(4m+4n-4)(3^a.2^b+3^b.2^a)+2(2^a.2^b+2^b.2^a).
\end{eqnarray*}
Hence, the theorem.                        \qed
 
\begin{cor}
From equation 1, we derived the following results,
\begin{eqnarray*}
(i)~M_1(H^f[m,n])\hspace{.75cm}&=&Z_{1,0}(H^f[m,n])=18mn+8m+8n+2,\\
(ii)~M_2(H^f[m,n])\hspace{.64cm}&=&\frac{1}{2} Z_{1,1}(H^f[m,n])=27mn+6m+6n-1,\\
(iii)~F(H^f[m,n])\hspace{.75cm}&=&Z_{2,0}(H^f[m,n])=54mn+16m+16n+2,\\
(iv)~ReZM(H^f[m,n])&=&Z_{2,1}(H^f[m,n])=162mn+12m+12n-10,\\
(v)~M^a(H^f[m,n])\hspace{.66cm}&=&Z_{a-1,0}(H^f[m,n])=6mn.3^{a-1}+2^{a+1}+2\\
&&+2^{a-1}(4m+4n-4),\\
(vi)~R_a(H^f[m,n])\hspace{.66cm}&=&\frac{1}{2} Z_{a,a}(H^f[m,n])=3^{2a}(3mn-2m-2n+1)+2^{2a+1}\\
&&+2^a.3^a(4m+4n-4)+2.3^a,\\
(vii)~SDD(H^f[m,n])&=&Z_{1,-1}(H^f[m,n])=2(3mn-2m-2n+1)\\
&&+\frac{13}{6}(4m+4n-4)+\frac{32}{3}.
\end{eqnarray*}
\end{cor}

Next, we consider the hexagonal lattice network with cylindrical boundary condition denoted as $H^c[m,n]$ and compute the general Zagreb index for this network. The edge partition of $H^c[m,n]$ based on degrees of the end vertices of the each edge is shown in the table 3. Note that the total number of edges of $H^c[m,n]$ is $(n+1)(3m+2)$ that is $|E(H^c[m,n])|=(n+1)(3m+2).$
\begin{table}[ht]
\caption{Edge partition of hexagonal lattice network with cylindrical boundary condition} 
\centering  
\begin{tabular}{|c|c|}
\hline
$(d(u),d(v)): ~uv\in E(H^c[m,n])$ & Total number of edges\\
\hline
(3,3) & $(n+1)(3m-2)$\\
\hline
(3,2) & $4(n+1)$\\
\hline
\end{tabular}
\label{table:nonlin} 
\end{table}

\begin{figure*}[ht]
\centering   
\captionsetup{singlelinecheck = true}
\begin{tikzpicture}[scale=.8,
every edge/.style = {draw=black, thick},
 vrtx/.style args = {#1/#2}{
      circle, draw, thick, fill=black,
      scale=.1, label=#1:#2}
                    ]
\node (A) [vrtx=/]   [label=left:{$b_1$}]   at (0,0) {};
\node (B) [vrtx=/]   [label=below left:{$a_1$}]  at (2,0) {};
\node (C) [vrtx=/]    [label=below left:{$a_2$}] at ( 4, 0) {};
\node (D) [vrtx=/]      at (6, 0) {};
\node (E) [vrtx=/]   [label=below left:{$a_{n}$}]  at ( 8, 0) {};
\node (F) [vrtx=/]    at ( 1, 0.5) {};
\node (G) [vrtx=/]     at ( 3, 0.5) {};
\node (H) [vrtx=/]     at ( 5, 0.5) {};
\node (I) [vrtx=/]    at ( 7, 0.5) {};
\node (J) [vrtx=/]   [label=below right:{$c_1$}] at ( 9, 0.5) {};
\node (K) [vrtx=/]    [label=left:{$b_2$}] at ( 1, 1.5) {};
\node (L) [vrtx=/]    at ( 3, 1.5) {};
\node (M) [vrtx=/]    at ( 5, 1.5) {};
\node (N) [vrtx=/]     at ( 7, 1.5) {};
\node (O) [vrtx=/]     at ( 9, 1.5) {};
\node (P) [vrtx=/]     at ( 2, 2) {};
\node (Q) [vrtx=/]     at ( 4, 2) {};
\node (R) [vrtx=/]     at ( 6, 2) {};
\node (S) [vrtx=/]      at ( 8, 2) {};
\node (T) [vrtx=/]    [label=below right:{$c_2$}] at (10,2) {};
\node (U) [vrtx=/]    [label=left:{$b_{3}$}] at (2,3) {};
\node (V) [vrtx=/]     at ( 4, 3) {};
\node (W) [vrtx=/]     at (6, 3) {};
\node (X) [vrtx=/]     at ( 8, 3) {};
\node (Y) [vrtx=/]     at ( 10, 3) {};
\node (Z) [vrtx=/]     at ( 3, 3.5) {};
\node (A*) [vrtx=/]     at ( 5, 3.5) {};
\node (B*) [vrtx=/]    at ( 7, 3.5) {};
\node (C*) [vrtx=/]    at ( 9, 3.5) {};
\node (D*) [vrtx=/]    at ( 11, 3.5) {};
\node (E*) [vrtx=/]     at ( 3, 4.5) {};
\node (F*) [vrtx=/]      at ( 5, 4.5) {};
\node (G*) [vrtx=/]     at ( 7, 4.5) {};
\node (H*) [vrtx=/]     at ( 9, 4.5) {};
\node (I*) [vrtx=/]     at ( 11, 4.5) {};
\node (J*) [vrtx=/]     at ( 4, 5) {};
\node (K*) [vrtx=/]       at ( 6, 5) {};
\node (L*) [vrtx=/]     at ( 8, 5) {};
\node (M*) [vrtx=/]     at ( 10, 5) {};
\node (N*) [vrtx=/]     [label=below right:{$c_{m}$}] at ( 12, 5) {};
\node (O*) [vrtx=/]     [label= left:{$b_{m+1}$}] at ( 4, 6) {};
\node (P*) [vrtx=/]     at ( 6, 6) {};
\node (Q*) [vrtx=/]     at ( 8, 6) {};
\node (R*) [vrtx=/]     at ( 10, 6) {};
\node (S*) [vrtx=/]     at ( 12, 6) {};
\node (T*) [vrtx=/]    [label=above left:{$a_{1}^*$}] at ( 5, 6.5) {};
\node (U*) [vrtx=/]     [label=above left:{$a_{2}^*$}] at ( 7, 6.5) {};
\node (V*) [vrtx=/]     at ( 9, 6.5) {};
\node (W*) [vrtx=/]     [label=above left:{$a_{n}^*$}] at ( 11, 6.5) {};
\node (X*) [vrtx=/]     [label=below right:{$c_{m+1}$}] at ( 13, 6.5) {};
\node (Y*) [vrtx=/]   at (5.9,1) {};
\node (Z*) [vrtx=/]   at (6,1) {};
\node (A**) [vrtx=/]   at (6.1,1) {};
\node (B**) [vrtx=/]   at (6.9,2.5) {};
\node (C**) [vrtx=/]   at (7,2.5) {};
\node (D**) [vrtx=/]   at (7.1,2.5) {};
\node (E**) [vrtx=/]   at (8.9,5.5) {};
\node (F**) [vrtx=/]   at (9,5.5) {};
\node (G**) [vrtx=/]   at (9.1,5.5) {};
\node (C**) [vrtx=/]   at (3.9,3.9) {};
\node (D**) [vrtx=/]   at (4,4) {};
\node (E**) [vrtx=/]   at (4.1,4.1) {};
\node (C**) [vrtx=/]   at (5.9,3.9) {};
\node (D**) [vrtx=/]   at (6,4) {};
\node (E**) [vrtx=/]   at (6.1,4.1) {};
\node (C**) [vrtx=/]   at (7.9,3.9) {};
\node (D**) [vrtx=/]   at (8,4) {};
\node (E**) [vrtx=/]   at (8.1,4.1) {};
\node (C**) [vrtx=/]   at (9.9,3.9) {};
\node (D**) [vrtx=/]   at (10,4) {};
\node (E**) [vrtx=/]   at (10.1,4.1) {};
\path   (A) edge (F)
        (B) edge (F)
        (B) edge (G)
        (G) edge (C)
        (C) edge (H)
        (H) edge (D)
        (D) edge (I)
        (I) edge (E)
        (E) edge (J)
        (F) edge (K)
        (G) edge (L)
        (H) edge (M)
        (I) edge (N)
        (J) edge (O)
        (K) edge (P)
        (P) edge (L)
        (L) edge (Q)
        (Q) edge (M)
        (M) edge (R)
        (R) edge (N)
        (N) edge (S)
        (S) edge (O)
        (O) edge (T)
        (P) edge (U)
        (Q) edge (V)
        (R) edge (W)
        (S) edge (X)
        (T) edge (Y)
        (U) edge (Z)
        (Z) edge (V)
        (V) edge (A*)
        (A*) edge (W)
        (W) edge (B*)
        (B*) edge (X)
        (X) edge (C*)
        (C*) edge (Y)
        (Y) edge (D*)
        (Z) edge (E*)
        (A*) edge (F*)
        (B*) edge (G*)
        (C*) edge (H*)
        (D*) edge (I*)
        (E*) edge (J*)
        (J*) edge (F*)
        (F*) edge (K*)
        (K*) edge (G*)
        (G*) edge (L*)
        (L*) edge (H*)
        (H*) edge (M*)
        (M*) edge (I*)
        (I*) edge (N*)
        (J*) edge (O*)
        (K*) edge (P*)
        (L*) edge (Q*)
        (M*) edge (R*)
        (N*) edge (S*)
        (O*) edge (T*)
        (T*) edge (P*)
        (P*) edge (U*)
        (U*) edge (Q*)
        (Q*) edge (V*)
        (V*) edge (R*)
        (R*) edge (W*)
        (W*) edge (S*)
        (S*) edge (X*);
      \path [ dotted]  (A) edge (J);
       \path [ dotted]  (K) edge (T);
       \path [ dotted]  (U) edge (D*);
     \path [ dotted]   (E*) edge (N*); 
      \path [ dotted]  (O*) edge (X*); 
\end{tikzpicture}
\caption*{~$H^c[m,n]$}
\captionsetup{singlelinecheck = true, justification=justified}
\caption{The example of hexagonal lattice with cylindrical boundary conditions.}
\end{figure*}

\begin{thm}
The general Zagreb index of $(H^c[m,n])$ is given by
\begin{eqnarray}
Z_{a,b}(H^c[m,n])&=&2(n+1)(3m-2).3^{(a+b)}+4(n+1)(3^a.2^b+3^b.2^a)
\end{eqnarray}  
\end{thm}
\textit{Proof.} By the definition of general Zagreb index, we get
\begin{eqnarray*}
Z_{a,b}(H^c[m,n])&=&\sum\limits_{uv\in E(H^c[m,n])}(d_{H^c[m,n]}(u)^ad_{H^c[m,n]}(v)^b+d_{H^c[m,n]}(u)^bd_{H^c[m,n]}(v)^a)
\end{eqnarray*}
\begin{eqnarray*}
&=&\sum\limits_{uv\in E_1(H^c[m,n])}(3^a.3^b+3^b.3^a)+\sum\limits_{uv\in E_2(H^c[m,n])}(3^a.2^b+3^b.2^a)\\
&=&|E_1(H^c[m,n])|(3^a.3^b+3^b.3^a)+|E_2(H^c[m,n])|(3^a.2^b+3^b.2^a)\\
&=&(n+1)(3m-2)(3^a.3^b+3^b.3^a)+4(n+1)(3^a.2^b+3^b.2^a).\\
\end{eqnarray*}
Hence, the result follows as in theorem 2.                        \qed
 
\begin{cor}
Using equation 2, we compute the following topological indices,
\begin{eqnarray*}
(i)~M_1(H^c[m,n])\hspace{.75cm}&=&Z_{1,0}(H^c[m,n])=2(n+1)(9m+4),\\
(ii)~M_2(H^c[m,n])\hspace{.64cm}&=&\frac{1}{2} Z_{1,1}(H^c[m,n])=3(n+1)(9m+2),\\
(iii)~F(H^c[m,n])\hspace{.75cm}&=&Z_{2,0}(H^c[m,n])=2(n+1)(27m+8),\\
(iv)~ReZM(H^c[m,n])&=&Z_{2,1}(H^c[m,n])=6(n+1)(27m+2),\\
(v)~M^a(H^c[m,n])\hspace{.66cm}&=&Z_{a-1,0}(H^c[m,n])=2(n+1)(3m-2).3^{a-1}\\
&&+4(n+1)(3^{a-1}+2^{a-1}),\\
(vi)~R_a(H^c[m,n])\hspace{.66cm}&=&\frac{1}{2} Z_{a,a}(H^c[m,n])=(n+1)(3m-2).3^{2a}+4(n+1).3^a.2^a,\\
(vii)~SDD(H^c[m,n])&=&Z_{1,-1}(H^c[m,n])=2(n+1)(3m+\frac{7}{3}).
\end{eqnarray*}
\end{cor}

Here, we compute the general Zagreb index of hexagonal lattice network with toroidal condition $H^t[m,n]$. In this network, it is clearly shown that the degree of all the vertices is the same and is equal to $3$ and the total number of edges in $H^t[m,n]$ is $3(m+1)(n+1)$ as shown in the table 4. 

\begin{table}[ht]
\caption{Edge partition of $H^t[m,n]$ network} 
\centering  
\begin{tabular}{|c|c|}
\hline
$(d(u),d(v)): ~uv\in E(H^t[m,n])$ & Total number of edges\\
\hline
(3,3) & 3(m+1)(n+1)\\
\hline
\end{tabular}
\label{table:nonlin} 
\end{table}

\begin{figure*}[ht]
\centering   
\captionsetup{singlelinecheck = true}
\begin{tikzpicture}[scale=.8,
every edge/.style = {draw=black, thick},
 vrtx/.style args = {#1/#2}{
      circle, draw, thick, fill=black,
      scale=.1, label=#1:#2}
                    ]
\node (A) [vrtx=/]   [label=left:{$b_1$}]   at (0,0) {};
\node (B) [vrtx=/]   [label=below left:{$a_1$}]  at (2,0) {};
\node (C) [vrtx=/]    [label=below left:{$a_2$}] at ( 4, 0) {};
\node (D) [vrtx=/]      at (6, 0) {};
\node (E) [vrtx=/]   [label=below left:{$a_{n}$}]  at ( 8, 0) {};
\node (F) [vrtx=/]    at ( 1, 0.5) {};
\node (G) [vrtx=/]     at ( 3, 0.5) {};
\node (H) [vrtx=/]     at ( 5, 0.5) {};
\node (I) [vrtx=/]    at ( 7, 0.5) {};
\node (J) [vrtx=/]   [label=below right:{$c_1$}] at ( 9, 0.5) {};
\node (K) [vrtx=/]    [label=left:{$b_2$}] at ( 1, 1.5) {};
\node (L) [vrtx=/]    at ( 3, 1.5) {};
\node (M) [vrtx=/]    at ( 5, 1.5) {};
\node (N) [vrtx=/]     at ( 7, 1.5) {};
\node (O) [vrtx=/]     at ( 9, 1.5) {};
\node (P) [vrtx=/]     at ( 2, 2) {};
\node (Q) [vrtx=/]     at ( 4, 2) {};
\node (R) [vrtx=/]     at ( 6, 2) {};
\node (S) [vrtx=/]      at ( 8, 2) {};
\node (T) [vrtx=/]    [label=below right:{$c_2$}] at (10,2) {};
\node (U) [vrtx=/]    [label=left:{$b_{3}$}] at (2,3) {};
\node (V) [vrtx=/]     at ( 4, 3) {};
\node (W) [vrtx=/]     at (6, 3) {};
\node (X) [vrtx=/]     at ( 8, 3) {};
\node (Y) [vrtx=/]     at ( 10, 3) {};
\node (Z) [vrtx=/]     at ( 3, 3.5) {};
\node (A*) [vrtx=/]     at ( 5, 3.5) {};
\node (B*) [vrtx=/]    at ( 7, 3.5) {};
\node (C*) [vrtx=/]    at ( 9, 3.5) {};
\node (D*) [vrtx=/]    at ( 11, 3.5) {};
\node (E*) [vrtx=/]     at ( 3, 4.5) {};
\node (F*) [vrtx=/]      at ( 5, 4.5) {};
\node (G*) [vrtx=/]     at ( 7, 4.5) {};
\node (H*) [vrtx=/]     at ( 9, 4.5) {};
\node (I*) [vrtx=/]     at ( 11, 4.5) {};
\node (J*) [vrtx=/]     at ( 4, 5) {};
\node (K*) [vrtx=/]       at ( 6, 5) {};
\node (L*) [vrtx=/]     at ( 8, 5) {};
\node (M*) [vrtx=/]     at ( 10, 5) {};
\node (N*) [vrtx=/]     [label=below right:{$c_{m}$}] at ( 12, 5) {};
\node (O*) [vrtx=/]     [label= left:{$b_{m+1}$}] at ( 4, 6) {};
\node (P*) [vrtx=/]     at ( 6, 6) {};
\node (Q*) [vrtx=/]     at ( 8, 6) {};
\node (R*) [vrtx=/]     at ( 10, 6) {};
\node (S*) [vrtx=/]     at ( 12, 6) {};
\node (T*) [vrtx=/]    [label=above left:{$a_{1}^*$}] at ( 5, 6.5) {};
\node (U*) [vrtx=/]     [label=above left:{$a_{2}^*$}] at ( 7, 6.5) {};
\node (V*) [vrtx=/]     at ( 9, 6.5) {};
\node (W*) [vrtx=/]     [label=above left:{$a_{n}^*$}] at ( 11, 6.5) {};
\node (X*) [vrtx=/]     [label=below right:{$c_{m+1}$}] at ( 13, 6.5) {};
\node (Y*) [vrtx=/]   at (5.9,1) {};
\node (Z*) [vrtx=/]   at (6,1) {};
\node (A**) [vrtx=/]   at (6.1,1) {};
\node (B**) [vrtx=/]   at (6.9,2.5) {};
\node (C**) [vrtx=/]   at (7,2.5) {};
\node (D**) [vrtx=/]   at (7.1,2.5) {};
\node (E**) [vrtx=/]   at (8.9,5.5) {};
\node (F**) [vrtx=/]   at (9,5.5) {};
\node (G**) [vrtx=/]   at (9.1,5.5) {};
\node (C**) [vrtx=/]   at (3.9,3.9) {};
\node (D**) [vrtx=/]   at (4,4) {};
\node (E**) [vrtx=/]   at (4.1,4.1) {};
\node (C**) [vrtx=/]   at (5.9,3.9) {};
\node (D**) [vrtx=/]   at (6,4) {};
\node (E**) [vrtx=/]   at (6.1,4.1) {};
\node (C**) [vrtx=/]   at (7.9,3.9) {};
\node (D**) [vrtx=/]   at (8,4) {};
\node (E**) [vrtx=/]   at (8.1,4.1) {};
\node (C**) [vrtx=/]   at (9.9,3.9) {};
\node (D**) [vrtx=/]   at (10,4) {};
\node (E**) [vrtx=/]   at (10.1,4.1) {};
\path   (A) edge (F)
        (B) edge (F)
        (B) edge (G)
        (G) edge (C)
        (C) edge (H)
        (H) edge (D)
        (D) edge (I)
        (I) edge (E)
        (E) edge (J)
        (F) edge (K)
        (G) edge (L)
        (H) edge (M)
        (I) edge (N)
        (J) edge (O)
        (K) edge (P)
        (P) edge (L)
        (L) edge (Q)
        (Q) edge (M)
        (M) edge (R)
        (R) edge (N)
        (N) edge (S)
        (S) edge (O)
        (O) edge (T)
        (P) edge (U)
        (Q) edge (V)
        (R) edge (W)
        (S) edge (X)
        (T) edge (Y)
        (U) edge (Z)
        (Z) edge (V)
        (V) edge (A*)
        (A*) edge (W)
        (W) edge (B*)
        (B*) edge (X)
        (X) edge (C*)
        (C*) edge (Y)
        (Y) edge (D*)
        (Z) edge (E*)
        (A*) edge (F*)
        (B*) edge (G*)
        (C*) edge (H*)
        (D*) edge (I*)
        (E*) edge (J*)
        (J*) edge (F*)
        (F*) edge (K*)
        (K*) edge (G*)
        (G*) edge (L*)
        (L*) edge (H*)
        (H*) edge (M*)
        (M*) edge (I*)
        (I*) edge (N*)
        (J*) edge (O*)
        (K*) edge (P*)
        (L*) edge (Q*)
        (M*) edge (R*)
        (N*) edge (S*)
        (O*) edge (T*)
        (T*) edge (P*)
        (P*) edge (U*)
        (U*) edge (Q*)
        (Q*) edge (V*)
        (V*) edge (R*)
        (R*) edge (W*)
        (W*) edge (S*)
        (S*) edge (X*);
      \path [ dotted]  (A) edge (J);
       \path [ dotted]  (K) edge (T);
       \path [ dotted]  (U) edge (D*);
     \path [ dotted]   (E*) edge (N*); 
      \path [ dotted]  (O*) edge (X*); 
       \path [ dotted]  (A) edge (T*); 
             \path [ dotted]  (B) edge (U*);
       \path [ dotted]  (C) edge (V*);
      \path [ dotted]  (D) edge (W*);
       \path [ dotted]  (E) edge (X*);
\end{tikzpicture}
\caption*{~$H^t[m,n]$}
\captionsetup{singlelinecheck = true, justification=justified}
\caption{The example of hexagonal lattice with toroidal boundary conditions.}
\end{figure*}

\begin{thm}
For the hexagonal lattice network with toroidal boundary condition $H^t[m,n]$, the general Zagreb index is given by
\begin{eqnarray}
 Z_{a,b}(H^t[m,n])&=&6(m+1)(n+1).3^{(a+b)}.
\end{eqnarray}  
\end{thm}
\textit{Proof.} We have from the definition of general Zagreb index, we get
\begin{eqnarray*}
Z_{a,b}(H^t[m,n])&=&\sum\limits_{uv\in E(H^t[m,n])}(d_{H^t[m,n]}(u)^ad_{H^t[m,n]}(v)^b+d_{H^t[m,n]}(u)^bd_{H^t[m,n]}(v)^a)\\
&=&\sum\limits_{uv\in E(H^t[m,n])}(3^a.3^b+3^b.3^a)
\end{eqnarray*}
\begin{eqnarray*}
&=&|E(H^t[m,n])|(3^a.3^b+3^b.3^a)\\
&=&3(m+1)(n+1)2.3^{(a+b)}.\\
\end{eqnarray*}
Hence, the theorem.                        \qed
 
\begin{cor}
Using equation 3, we derived the following results,
\begin{eqnarray*}
(i)~M_1(H^t[m,n])\hspace{.75cm}&=&Z_{1,0}(H^t[m,n])=18(m+1)(n+1),\\
(ii)~M_2(H^t[m,n])\hspace{.64cm}&=&\frac{1}{2} Z_{1,1}(H^t[m,n])=27(m+1)(n+1),\\
(iii)~F(H^t[m,n])\hspace{.75cm}&=&Z_{2,0}(H^t[m,n])=54(m+1)(n+1),\\
(iv)~ReZM(H^t[m,n])&=&Z_{2,1}(H^t[m,n])=162(m+1)(n+1),\\
(v)~M^a(H^t[m,n])\hspace{.66cm}&=&Z_{a-1,0}(H^t[m,n])=6(m+1)(n+1).3^{a-1},\\
(vi)~R_a(H^t[m,n])\hspace{.66cm}&=&\frac{1}{2} Z_{a,a}(H^t[m,n])=(m+1)(n+1).3^{2a+1},\\
(vii)~SDD(H^t[m,n])&=&Z_{1,-1}(H^t[m,n])=6(m+1)(n+1).
\end{eqnarray*}
\end{cor}
Now, we consider the triangular lattice network with free boundary condition $T^f[m,n]$, such that $|E(T^f[m,n])|=m(n-1)+n(m-1)+(n-1)(m-1)$. The edge partition of $T^f[m,n]$ network with respect to the degree of the end vertices of each edge is shown in the table 5, as follows: 
\begin{table}[ht]
\caption{Edge partition of $T^f[m,n]$ network} 
\centering  
\begin{tabular}{|c|c|}
\hline
$(d(u),d(v)): ~uv\in E(T^f[m,n])$ & Total number of edges\\
\hline
(2,4) & $4$\\
\hline
(3,4) & $4$\\
\hline
(3,6) & $2$\\
\hline
(4,4) & $2(m+n-5)$\\
\hline
(4,6) & $4(m+n-5)$\\
\hline
(6,6) & $3(mn+7)-8(m+n)$\\
\hline
\end{tabular}
\label{table:nonlin} 
\end{table}

\begin{figure*}[h]
\centering   
\captionsetup{singlelinecheck = true}
\begin{tikzpicture}[scale=.8,
every edge/.style = {draw=black, thick},
 vrtx/.style args = {#1/#2}{
      circle, draw, thick, fill=black,
      scale=.2, label=#1:#2}
                    ]
\node (A) [vrtx=/]   [label=below right:{$a_1$}] [label=left:{$b_1$}]  at (1,2) {};
\node (B) [vrtx=/]   [label=below right:{$a_2$}]  at (3,2) {};
\node (C) [vrtx=/]     at ( 5, 2) {};
\node (D) [vrtx=/]     at (7, 2) {};
\node (E) [vrtx=/]   [label=below right:{$a_{n-1}$}]  at ( 9, 2) {};
\node (F) [vrtx=/]    [label=below right:{$a_n$}] [label=above left:{$b_1^*$}]at ( 11, 2) {};
\node (G) [vrtx=/]    [label=above left:{$b_2$}] at ( 1, 4) {};
\node (H) [vrtx=/]     at ( 3, 4) {};
\node (I) [vrtx=/]    at ( 5, 4) {};
\node (J) [vrtx=/]    at ( 7, 4) {};
\node (K) [vrtx=/]    at ( 9, 4) {};
\node (L) [vrtx=/]    [label=above left:{$b_2^*$}] at ( 11, 4) {};
\node (M) [vrtx=/]     at ( 1, 6) {};
\node (N) [vrtx=/]     at ( 3, 6) {};
\node (O) [vrtx=/]     at ( 5, 6) {};
\node (P) [vrtx=/]     at ( 7, 6) {};
\node (Q) [vrtx=/]     at ( 9, 6) {};
\node (R) [vrtx=/]      at ( 11, 6) {};
\node (S) [vrtx=/]     at (1,8) {};
\node (T) [vrtx=/]     at (3,8) {};
\node (U) [vrtx=/]     at ( 5, 8) {};
\node (V) [vrtx=/]     at (7, 8) {};
\node (W) [vrtx=/]     at ( 9, 8) {};
\node (X) [vrtx=/]     at ( 11, 8) {};
\node (Y) [vrtx=/]    [label=above left:{$b_{m-1}$}] at ( 1, 10) {};
\node (Z) [vrtx=/]     at ( 3, 10) {};
\node (A*) [vrtx=/]    at ( 5, 10) {};
\node (B*) [vrtx=/]    at ( 7, 10) {};
\node (C*) [vrtx=/]    at ( 9, 10) {};
\node (D*) [vrtx=/]    [label=above left:{$b_{m-1}^*$}] at ( 11, 10) {};
\node (E*) [vrtx=/]    [label=above left:{$b_m$}] [label=below right :{$a_1^*$}] at ( 1, 12) {};
\node (F*) [vrtx=/]   [label=below right :{$a_2^*$}]  at ( 3, 12) {};
\node (G*) [vrtx=/]     at ( 5, 12) {};
\node (H*) [vrtx=/]     at ( 7, 12) {};
\node (I*) [vrtx=/]    [label=below right :{$a_{n-1}^*$}] at ( 9, 12) {};
\node (J*) [vrtx=/]    [label=above left:{$b_m^*$}] [label=below right :{$a_n^*$}]   at ( 11, 12) {};
\node (K*) [vrtx=/]     at ( 5.7, 1.8) {};
\node (L*) [vrtx=/]     at ( 5.9, 1.8) {};
\node (M*) [vrtx=/]     at ( 6.1, 1.8) {};
\node (N*) [vrtx=/]     at ( 7.7, 1.8) {};
\node (O*) [vrtx=/]     at ( 7.9, 1.8) {};
\node (P*) [vrtx=/]     at ( 8.1, 1.8) {};
\node (Q*) [vrtx=/]     at ( 5.7, 11.8) {};
\node (R*) [vrtx=/]     at ( 5.9, 11.8) {};
\node (S*) [vrtx=/]     at ( 6.1, 11.8) {};
\node (T*) [vrtx=/]     at ( 7.7, 11.8) {};
\node (U*) [vrtx=/]     at ( 7.9, 11.8) {};
\node (V*) [vrtx=/]     at ( 8.1, 11.8) {};

\path   (A) edge (F)
        (G) edge (L)
        (M) edge (R)
        (S) edge (X)
        (Y) edge (D*)
        (E*) edge (J*)
        (A) edge (E*)
        (B) edge (F*)
        (C) edge (G*)
        (D) edge (H*)
        (E) edge (I*)
        (F) edge (J*)
        (E) edge (L)
        (D) edge (R)
        (C) edge (X)
        (B) edge (D*)
        (A) edge (J*)
        (G) edge (I*)
        (M) edge (H*)
        (S) edge (G*)
        (Y) edge (F*);

\end{tikzpicture}
\caption*{~$T^f[m,n]$}
\captionsetup{singlelinecheck = true, justification=justified}
\caption{The example of triangular lattice with free boundary conditions.}
\end{figure*}

\begin{thm}The general Zagreb index of  $T^f[m,n]$ is given by
\begin{eqnarray}
\nonumber Z_{a,b}(T^f[m,n])&=&4(2^{a+2b}+2^{2a+b})+2[3(mn+7)-8(m+n)].6^{a+b}\\
\nonumber &&+2(3^a.6^b+3^b.6^a)+(m+n-5)2^{2(a+b+1)}\\
&&+4(m+n-5)(4^a.6^b+4^b.6^a)+4(3^a.4^b+3^b.4^a).
\end{eqnarray}  
\end{thm}
\textit{Proof.} We have from the definition of general Zagreb index, we get
\begin{eqnarray*}
Z_{a,b}(T^f[m,n])&=&\sum\limits_{uv\in E(T^f[m,n])}[d_{T^f[m,n]}(u)^ad_{T^f[m,n]}(v)^b+d_{T^f[m,n]}(u)^bd_{T^f[m,n]}(v)^a]\\
&=&\sum\limits_{uv\in E_1(T^f[m,n])}(2^a.4^b+2^b.4^a)+\sum\limits_{uv\in E_2(T^f[m,n])}(3^a.4^b+3^b.4^a)\\
&&+\sum\limits_{uv\in E_3(T^f[m,n])}(3^a.6^b+3^b.6^a)+\sum\limits_{uv\in E_4(T^f[m,n])}(4^a.4^b+4^b.4^a)\\
&&+\sum\limits_{uv\in E_5(T^f[m,n])}(4^a.6^b+4^b.6^a)+\sum\limits_{uv\in E_6(T^f[m,n])}(6^a.6^b+6^b.6^a)\\
&=&|E_1(T^f[m,n])|(2^a.4^b+2^b.4^a)+|E_2(T^f[m,n])|(3^a.4^b+3^b.4^a)\\
&&+|E_3(T^f[m,n])|(3^a.6^b+3^b.6^a)+|E_4(T^f[m,n])|(4^a.4^b+4^b.4^a)\\
&&+|E_5(T^f[m,n])|(4^a.6^b+4^b.6^a)+|E_6(T^f[m,n])|(6^a.6^b+6^b.6^a)\\
&=&4(2^a.4^b+2^b.4^a)+[3(mn+7)-8(m+n)](6^a.6^b+6^b.6^a)
\end{eqnarray*}
\begin{eqnarray*}
&&+2(3^a.6^b+3^b.6^a)+2(m+n-5)(4^a.4^b+4^b.4^a)\\
&&+4(m+n-5)(4^a.6^b+4^b.6^a)+4(3^a.4^b+3^b.4^a).
\end{eqnarray*}
Which is the required result as in theorem 4.                        \qed
 
\begin{cor}
Using equation 4, we derived the following results,
\begin{eqnarray*}
(i)~M_1(T^f[m,n])\hspace{.75cm}&=&Z_{1,0}(T^f[m,n])=12[3(mn+7)-8(m+n)]\\
&&+56(m+n-5)+70,\\
(ii)~M_2(T^f[m,n])\hspace{.64cm}&=&\frac{1}{2} Z_{1,1}(T^f[m,n])=36[3(mn+7)-8(m+n)]\\
&&+128(m+n-5)+116,\\
(iii)~F(T^f[m,n])\hspace{.75cm}&=&Z_{2,0}(T^f[m,n])=72[3(mn+7)-8(m+n)]\\
&&+272(m+n-5)+270,\\
(iv)~ReZM(T^f[m,n])&=&Z_{2,1}(T^f[m,n])=432[3(mn+7)-8(m+n)]\\
&&+1216(m+n-5)+852,
\end{eqnarray*}
\begin{eqnarray*}
(v)~M^a(T^f[m,n])\hspace{.66cm}&=&Z_{a-1,0}(T^f[m,n])=2[3(mn+7)-8(m+n)].6^{a-1}\\
&&+2(3^{a-1}+6^{a-1})+4(3^{a-1}+4^{a-1})+4(2^{a-1}+2^{2(a-1)})\\
&&+(m+n-5)[2^{2a}+4(4^{a-1}+6^{a-1})],\\
(vi)~R_a(T^f[m,n])\hspace{.66cm}&=&\frac{1}{2} Z_{a,a}(T^f[m,n])=[3(mn+7)-8(m+n)].6^{2a}+4.2^{3a}\\
&&+(m+n-5).2^{4a+1}+2.3^a.6^a\\
&&+4(m+n-5).4^a.6^a+4.3^a.4^a,\\
(vii)~SDD(T^f[m,n])&=&Z_{1,-1}(T^f[m,n])=2[3(mn+7)-8(m+n)]\\
&&+\frac{38}{3}(m+n-5)+\frac{70}{3}.
\end{eqnarray*}
\end{cor}

Next, we compute the general Zagreb index of triangular lattice with cylindrical boundary condition $T^c[m,n]$. Based on degree of the end vertices of each edge the edge partitioned of $T^c[m,n]$ is shown in table 6. Note that $|E(T^c[m,n])|=3mn-2n$.
\begin{table}[ht]
\caption{Edge partition of $T^c[m,n]$} 
\centering  
\begin{tabular}{|c|c|}
\hline
$(d(u),d(v)): ~uv\in E(T^c[m,n])$ & Total number of edges\\
\hline
(4,4) & $2n$\\
\hline
(4,6) & $4n$\\
\hline
(6,6) & $3mn-8n$\\
\hline
\end{tabular}
\label{table:nonlin} 
\end{table}

\begin{figure*}[ht]
\centering   
\captionsetup{singlelinecheck = true}
\begin{tikzpicture}[scale=.8,
every edge/.style = {draw=black, thick},
 vrtx/.style args = {#1/#2}{
      circle, draw, thick, fill=black,
      scale=.2, label=#1:#2}
                    ]
\node (A) [vrtx=/]   [label=below right:{$a_1$}] [label=left:{$b_1$}]  at (1,2) {};
\node (B) [vrtx=/]   [label=below right:{$a_2$}]  at (3,2) {};
\node (C) [vrtx=/]     at ( 5, 2) {};
\node (D) [vrtx=/]     at (7, 2) {};
\node (E) [vrtx=/]   [label=below right:{$a_{n-1}$}]  at ( 9, 2) {};
\node (F) [vrtx=/]    [label=below right:{$a_n$}] [label=above left:{$b_1^*$}]at ( 11, 2) {};
\node (G) [vrtx=/]    [label=above left:{$b_2$}] at ( 1, 4) {};
\node (H) [vrtx=/]     at ( 3, 4) {};
\node (I) [vrtx=/]    at ( 5, 4) {};
\node (J) [vrtx=/]    at ( 7, 4) {};
\node (K) [vrtx=/]    at ( 9, 4) {};
\node (L) [vrtx=/]    [label=above left:{$b_2^*$}] at ( 11, 4) {};
\node (M) [vrtx=/]     at ( 1, 6) {};
\node (N) [vrtx=/]     at ( 3, 6) {};
\node (O) [vrtx=/]     at ( 5, 6) {};
\node (P) [vrtx=/]     at ( 7, 6) {};
\node (Q) [vrtx=/]     at ( 9, 6) {};
\node (R) [vrtx=/]      at ( 11, 6) {};
\node (S) [vrtx=/]     at (1,8) {};
\node (T) [vrtx=/]     at (3,8) {};
\node (U) [vrtx=/]     at ( 5, 8) {};
\node (V) [vrtx=/]     at (7, 8) {};
\node (W) [vrtx=/]     at ( 9, 8) {};
\node (X) [vrtx=/]     at ( 11, 8) {};
\node (Y) [vrtx=/]    [label=above left:{$b_{m-1}$}] at ( 1, 10) {};
\node (Z) [vrtx=/]     at ( 3, 10) {};
\node (A*) [vrtx=/]    at ( 5, 10) {};
\node (B*) [vrtx=/]    at ( 7, 10) {};
\node (C*) [vrtx=/]    at ( 9, 10) {};
\node (D*) [vrtx=/]    [label=above left:{$b_{m-1}^*$}] at ( 11, 10) {};
\node (E*) [vrtx=/]    [label=above left:{$b_m$}] [label=below right :{$a_1^*$}] at ( 1, 12) {};
\node (F*) [vrtx=/]   [label=below right :{$a_2^*$}]  at ( 3, 12) {};
\node (G*) [vrtx=/]     at ( 5, 12) {};
\node (H*) [vrtx=/]     at ( 7, 12) {};
\node (I*) [vrtx=/]    [label=below right :{$a_{n-1}^*$}] at ( 9, 12) {};
\node (J*) [vrtx=/]    [label=above left:{$b_m^*$}] [label=below right :{$a_n^*$}]   at ( 11, 12) {};
\node (K*) [vrtx=/]     at ( 5.7, 1.8) {};
\node (L*) [vrtx=/]     at ( 5.9, 1.8) {};
\node (M*) [vrtx=/]     at ( 6.1, 1.8) {};
\node (N*) [vrtx=/]     at ( 7.7, 1.8) {};
\node (O*) [vrtx=/]     at ( 7.9, 1.8) {};
\node (P*) [vrtx=/]     at ( 8.1, 1.8) {};
\node (Q*) [vrtx=/]     at ( 5.7, 11.8) {};
\node (R*) [vrtx=/]     at ( 5.9, 11.8) {};
\node (S*) [vrtx=/]     at ( 6.1, 11.8) {};
\node (T*) [vrtx=/]     at ( 7.7, 11.8) {};
\node (U*) [vrtx=/]     at ( 7.9, 11.8) {};
\node (V*) [vrtx=/]     at ( 8.1, 11.8) {};

\path   (A) edge (F)
        (G) edge (L)
        (M) edge (R)
        (S) edge (X)
        (Y) edge (D*)
        (E*) edge (J*)
        (A) edge (E*)
        (B) edge (F*)
        (C) edge (G*)
        (D) edge (H*)
        (E) edge (I*)
        (F) edge (J*)
        (E) edge (L)
        (D) edge (R)
        (C) edge (X)
        (B) edge (D*)
        (A) edge (J*)
        (G) edge (I*)
        (M) edge (H*)
        (S) edge (G*)
        (Y) edge (F*);
       \path [ dotted] [in=160, out=20] (E*) edge (J*);
       \path [ dotted] [in=160, out=20] (Y) edge (D*);
       \path [ dotted] [in=160, out=20] (S) edge (X); 
       \path [ dotted] [in=160, out=20] (M) edge (R); 
       \path [ dotted] [in=160, out=20] (G) edge (L); 
       \path [ dotted] [in=160, out=20] (A) edge (F);
      \path [ dotted] (G) edge (F);
       \path [ dotted] (M) edge (L);
        \path [ dotted]  (S) edge (R);
       \path [ dotted]  (Y) edge (X);
        \path [ dotted]  (E*) edge (D*);
\end{tikzpicture}
\caption*{~$T^c[m,n]$}
\captionsetup{singlelinecheck = true, justification=justified}
\caption{The example of triangular lattice with cylindrical boundary conditions.}
\end{figure*}
\begin{thm} The general Zagreb index of $T^c[m,n]$ is given by
\begin{eqnarray}
Z_{a,b}(T^c[m,n])&=&n.2^{2(a+b+1)}+4n(4^a.6^b+4^b.6^a)+2(3mn-8n)6^{(a+b)}.
\end{eqnarray}  
\end{thm}
\textit{Proof.} By the definition of general Zagreb index, we get
\begin{eqnarray*}
Z_{a,b}(T^c[m,n])&=&\sum\limits_{uv\in E(T^c[m,n])}[d_{T^c[m,n]}(u)^ad_{T^c[m,n]}(v)^b+d_{T^c[m,n]}(u)^bd_{T^c[m,n]}(v)^a]\\
&=&\sum\limits_{uv\in E_1(T^c[m,n])}(4^a.4^b+4^b.4^a)+\sum\limits_{uv\in E_2(T^c[m,n])}(4^a.6^b+4^b.6^a)\\
&&+\sum\limits_{uv\in E_3(T^c[m,n])}(6^a.6^b+6^b.6^a)\\
&=&|E_1(T^c[m,n])|(4^a.4^b+4^b.4^a)+|E_2(T^c[m,n])|(4^a.6^b+4^b.6^a)\\
&&+|E_3(T^c[m,n])|(6^a.6^b+6^b.6^a)\\
&=&2n(4^a.4^b+4^b.4^a)+4n(4^a.6^b+4^b.6^a)\\
&&+(3mn-8n)(6^a.6^b+6^b.6^a).
\end{eqnarray*}
Hence the theorem.                        \qed
 
\begin{cor}
Using equation 5, we derived the following results,
\begin{eqnarray*}
(i)~M_1(T^c[m,n])\hspace{.75cm}&=&Z_{1,0}(T^c[m,n])=56n+12(3mn-8n),\\
(ii)~M_2(T^c[m,n])\hspace{.64cm}&=&\frac{1}{2} Z_{1,1}(T^c[m,n])=128n+36(3mn-8n),\\
(iii)~F(T^c[m,n])\hspace{.75cm}&=&Z_{2,0}(T^c[m,n])=272n+72(3mn-8n),\\
(iv)~ReZM(T^c[m,n])&=&Z_{2,1}(T^c[m,n])=1216n+432(3mn-8n),\\
(v)~M^a(T^c[m,n])\hspace{.66cm}&=&Z_{a-1,0}(T^c[m,n])=n.2^{2a}+4n(4^{a-1}+6^{a-1})\\
&&+2(3mn-8n).6^{a-1},
\end{eqnarray*}
\begin{eqnarray*}
(vi)~R_a(T^c[m,n])\hspace{.66cm}&=&\frac{1}{2} Z_{a,a}(T^c[m,n])=n.2^{4a+1}+4n.4^a.6^a\\
&&+(3mn-8n).6^{2a},\\
(vii)~SDD(T^c[m,n])&=&Z_{1,-1}(T^c[m,n])=\frac{38}{3}n+2(3mn-8n).
\end{eqnarray*}
\end{cor}
 
Finally, we consider the triangular lattice with toroidal boundary condition $T^t[m,n]$. The edge partition of this network based on the degree of the end vertices of each edge is shown in table 7. 

\begin{table}[ht]
\caption{Edge partition of $T^t[m,n]$ network} 
\centering  
\begin{tabular}{|c|c|}
\hline
$(d(u),d(v)): ~uv\in E(T^t[m,n])$ & Total number of edges\\
\hline
(6,6) & 3mn\\
\hline
\end{tabular}
\label{table:nonlin} 
\end{table}

\begin{figure*}[ht]
\centering   
\captionsetup{singlelinecheck = true}
\begin{tikzpicture}[scale=.8,
every edge/.style = {draw=black, thick},
 vrtx/.style args = {#1/#2}{
      circle, draw, thick, fill=black,
      scale=.2, label=#1:#2}
                    ]
\node (A) [vrtx=/]   [label=below right:{$a_1$}] [label=left:{$b_1$}]  at (1,2) {};
\node (B) [vrtx=/]   [label=below right:{$a_2$}]  at (3,2) {};
\node (C) [vrtx=/]     at ( 5, 2) {};
\node (D) [vrtx=/]     at (7, 2) {};
\node (E) [vrtx=/]   [label=below right:{$a_{n-1}$}]  at ( 9, 2) {};
\node (F) [vrtx=/]    [label=below right:{$a_n$}] [label=above left:{$b_1^*$}]at ( 11, 2) {};
\node (G) [vrtx=/]    [label=above left:{$b_2$}] at ( 1, 4) {};
\node (H) [vrtx=/]     at ( 3, 4) {};
\node (I) [vrtx=/]    at ( 5, 4) {};
\node (J) [vrtx=/]    at ( 7, 4) {};
\node (K) [vrtx=/]    at ( 9, 4) {};
\node (L) [vrtx=/]    [label=above left:{$b_2^*$}] at ( 11, 4) {};
\node (M) [vrtx=/]     at ( 1, 6) {};
\node (N) [vrtx=/]     at ( 3, 6) {};
\node (O) [vrtx=/]     at ( 5, 6) {};
\node (P) [vrtx=/]     at ( 7, 6) {};
\node (Q) [vrtx=/]     at ( 9, 6) {};
\node (R) [vrtx=/]      at ( 11, 6) {};
\node (S) [vrtx=/]     at (1,8) {};
\node (T) [vrtx=/]     at (3,8) {};
\node (U) [vrtx=/]     at ( 5, 8) {};
\node (V) [vrtx=/]     at (7, 8) {};
\node (W) [vrtx=/]     at ( 9, 8) {};
\node (X) [vrtx=/]     at ( 11, 8) {};
\node (Y) [vrtx=/]    [label=above left:{$b_{m-1}$}] at ( 1, 10) {};
\node (Z) [vrtx=/]     at ( 3, 10) {};
\node (A*) [vrtx=/]    at ( 5, 10) {};
\node (B*) [vrtx=/]    at ( 7, 10) {};
\node (C*) [vrtx=/]    at ( 9, 10) {};
\node (D*) [vrtx=/]    [label=above left:{$b_{m-1}^*$}] at ( 11, 10) {};
\node (E*) [vrtx=/]    [label=above left:{$b_m$}] [label=below right :{$a_1^*$}] at ( 1, 12) {};
\node (F*) [vrtx=/]   [label=below right :{$a_2^*$}]  at ( 3, 12) {};
\node (G*) [vrtx=/]     at ( 5, 12) {};
\node (H*) [vrtx=/]     at ( 7, 12) {};
\node (I*) [vrtx=/]    [label=below right :{$a_{n-1}^*$}] at ( 9, 12) {};
\node (J*) [vrtx=/]    [label=above left:{$b_m^*$}] [label=below right :{$a_n^*$}]   at ( 11, 12) {};
\node (K*) [vrtx=/]     at ( 5.7, 1.8) {};
\node (L*) [vrtx=/]     at ( 5.9, 1.8) {};
\node (M*) [vrtx=/]     at ( 6.1, 1.8) {};
\node (N*) [vrtx=/]     at ( 7.7, 1.8) {};
\node (O*) [vrtx=/]     at ( 7.9, 1.8) {};
\node (P*) [vrtx=/]     at ( 8.1, 1.8) {};
\node (Q*) [vrtx=/]     at ( 5.7, 11.8) {};
\node (R*) [vrtx=/]     at ( 5.9, 11.8) {};
\node (S*) [vrtx=/]     at ( 6.1, 11.8) {};
\node (T*) [vrtx=/]     at ( 7.7, 11.8) {};
\node (U*) [vrtx=/]     at ( 7.9, 11.8) {};
\node (V*) [vrtx=/]     at ( 8.1, 11.8) {};

\path   (A) edge (F)
        (G) edge (L)
        (M) edge (R)
        (S) edge (X)
        (Y) edge (D*)
        (E*) edge (J*)
        (A) edge (E*)
        (B) edge (F*)
        (C) edge (G*)
        (D) edge (H*)
        (E) edge (I*)
        (F) edge (J*)
        (E) edge (L)
        (D) edge (R)
        (C) edge (X)
        (B) edge (D*)
        (A) edge (J*)
        (G) edge (I*)
        (M) edge (H*)
        (S) edge (G*)
        (Y) edge (F*);
       \path [ dotted] [in=210, out=55] (A) edge (J*);
       \path [ dotted] [in=160, out=20] (E*) edge (J*);
       \path [ dotted] [in=160, out=20] (Y) edge (D*);
       \path [ dotted] [in=160, out=20] (S) edge (X); 
       \path [ dotted] [in=160, out=20] (M) edge (R); 
       \path [ dotted] [in=160, out=20] (G) edge (L); 
       \path [ dotted] [in=160, out=20] (A) edge (F);
       \path [ dotted] [in=220, out=125] (A) edge (E*);
        \path [ dotted] [in=250, out=110] (B) edge (F*);
        \path [ dotted] [in=250, out=110] (C) edge (G*);
        \path [ dotted] [in=250, out=110] (D) edge (H*);
        \path [ dotted] [in=250, out=110] (E) edge (I*);
        \path [ dotted] [in=250, out=110] (F) edge (J*);
        \path [ dotted] (G) edge (F);
       \path [ dotted] (M) edge (L);
        \path [ dotted]  (S) edge (R);
       \path [ dotted]  (Y) edge (X);
        \path [ dotted]  (E*) edge (D*);
        \path [ dotted] (B) edge (E*);
       \path [ dotted] (C) edge (F*);
        \path [ dotted]  (D) edge (G*);
       \path [ dotted]  (E) edge (H*);
        \path [ dotted]  (F) edge (I*);
\end{tikzpicture}
\caption*{~$T^t[m,n]$}
\captionsetup{singlelinecheck = true, justification=justified}
\caption{The example of triangular lattice with toroidal boundary conditions.}
\end{figure*}
\begin{thm}
For the triangular lattice network with toroidal boundary condition $T^t[m,n]$, the general Zagreb index is given by
\begin{eqnarray}
 Z_{a,b}(T^t[m,n])&=&6mn.6^{(a+b)}.
\end{eqnarray}  
\end{thm}
\textit{Proof.} We have from the definition of general Zagreb index, we get
\begin{eqnarray*}
Z_{a,b}(T^t[m,n])&=&\sum\limits_{uv\in E(T^t[m,n])}(d_{T^t[m,n]}(u)^ad_{T^t[m,n]}(v)^b+d_{T^t[m,n]}(u)^bd_{T^t[m,n]}(v)^a)\\
&=&\sum\limits_{uv\in E(T^t[m,n])}(6^a.6^b+6^b.6^a)\\
&=&|E(T^t[m,n])|(6^a.6^b+6^b.6^a)\\
&=&3mn(6^a.6^b+6^b.6^a).\\
\end{eqnarray*}
Hence, the theorem.                        \qed
 
\begin{cor}
Using equation 7, we derived the following results,
\begin{eqnarray*}
(i)~M_1(T^t[m,n])\hspace{.75cm}&=&Z_{1,0}(T^t[m,n])=36mn,\\
(ii)~M_2(T^t[m,n])\hspace{.64cm}&=&\frac{1}{2} Z_{1,1}(T^t[m,n])=18mn,\\
(iii)~F(H^t[m,n])\hspace{.75cm}&=&Z_{2,0}(H^t[m,n])=216mn,\\
(iv)~ReZM(T^t[m,n])&=&Z_{2,1}(T^t[m,n])=1296mn,\\
(v)~M^a(T^t[m,n])\hspace{.66cm}&=&Z_{a-1,0}(T^t[m,n])=6mn.6^{a-1},\\
(vi)~R_a(T^t[m,n])\hspace{.66cm}&=&\frac{1}{2} Z_{a,a}(T^t[m,n])=3mn.6^{2a},\\
(vii)~SDD(T^t[m,n])&=&Z_{1,-1}(T^t[m,n])=6mn.
\end{eqnarray*}
\end{cor}

\section{Conclusions}
 In this paper, the general Zagreb index or $(a,b)$-Zagreb index of hexagonal and triangular lattice network with three boundary conditions namely toroidal, cylindrical and free boundary conditions are obtained. Based on this general Zagreb index or $(a,b)$-Zagreb index of these network, we compute some other vertex degree based topological indices such as Zagreb indices, forgotten topological indices, redefined Zagreb index, general first Zagreb index, general Randi\'{c} index and symmetric division deg index for some particular values of $a$ and $b$. For future study, we want to compute some other topological indices of these networks.  

\section*{Acknowledgement} The first author would like to express sincere gratitude to CSIR, HRDG, New Delhi, India for their financial support under the grants no. 09/973(0016)/2017-EMR-I.



\begin{thebibliography}{99}

\bibitem{gutm72}I. Gutman ,  N. Trinajesti\'{c} ,          \textit{Graph theory and molecular orbitals total $\pi$-electron energy of alternant hydrocarbons},         {Chem. Phys. Lett.},      \textbf{17}, 535-538. 1972.

\bibitem{gut18}I. Gutman,  E. Milovanovic, and I. Milovanovic,    \textit{Beyond the Zagreb indices},          {AKCE Int. J. Graphs Comb.},       Doi: 10.1016/j.akcej.2018.05.002.  2018.

\bibitem{ndep17}P. Sarkar,  N. De, and  A. Pal,       \textit{The Zagreb indices of graphs based on new operations  related to the join of graphs},      {J. Int. Math. Virtual Inst.},     \textbf{7},    181-209.  2017.

\bibitem{sidd16}M.K. Siddiqui,  M. Imran,  and A. Ahmad,    \textit{On Zagreb indices, Zagreb polynomials of some nanostar dendrimers},        {Appl. Math. Comput.},    \textbf{280},  132-139.  2016.

\bibitem{furt15}B. Furtul,  I. Gutman,                 \textit{A forgotten topological index},            {J. Math. Chem.},               \textbf{53}, 1184-1190.   2015.

\bibitem{den17}De N.,       {F-index of total transformation graphs},             \textit{Discrete Math. Algorithm. Appl.},              \textbf{9},(3),(2017).

\bibitem{ram16}H.S. Ramane,  R.B. Jummannaver,     \textit{Note on forgotten topological index of chemical structures in drugs},        {Appl. Math. Nonlinear Sci.},     \textbf{1}, (2), 369-374.  2016.

\bibitem{akh17}S. Akhter,  M. Imran,     \textit{Computing the forgotten topological index of four operations on graphs},      {AKCE Int. J. Graphs Comb.},     \textbf{14}, (1),   70-79.   2017.

\bibitem{thay18}P.J.N. Thayamathy,  P. Elango, and M. koneswaran,   \textit{Forgotten topological index of some nano structures},       {Int. J. Current Mult. Stud.},    \textbf{4}, (7(A)),   909-914.   2018.

\bibitem{ranj13}P.S. Ranjini,    V. Lokesha, and  A. Usha,         \textit{Relation between phenylene and hexagonal squeeze using harmonic index},                {Int. J. Graph Theory},           \textbf{1}, 116-121.  2013.

\bibitem{bzha16}B. Zhao,  J. Gan, and  H. Wu,     \textit{Redefined Zagreb indices of some nano structures},      {Appl. Math. Nonlinear Sci.},     \textbf{1}, (1),   291-300.  2016.

\bibitem{gao16}W. Gao, M.R. Farahani, M.K. Jamil, and  M.K. Siddiqui,   \textit{The redefined first, second and third Zagreb indices of titania nanotubes $Tio_2[m,n]$},    {The Open Biotech. J.},     \textbf{10},  272-277.  2016. 

\bibitem{kum17}R.P. Kumar,   D.S. Nandappa, and   M.R.R. Kanna,      \textit{Redefined zagreb, Randi\'{c}, Harmonic, GA indices of graphene},     {Int. J. Math. Anal.},      \textbf{11}, (10),   493-502.  2017.

\bibitem{li05}X. Li,  J. Zheng,       \textit{A unified approach to the extremal trees for different indices},  {MATCH Commun. Math. Comput. Chem.},     \textbf{54}, 195-208.   2005.

\bibitem{liu10}M. Liu, B. Liu,     \textit{Some properties of the first general Zagreb index},      {Aus. J. Comb.},     \textbf{47}, 285-294.   2010.

\bibitem{bed18}L. Bedratyuk, O. Savenko,    \textit{The star sequence and the general first Zagreb index},    {MATCH Commun. Math. Comput. Chem.},    \textbf{79},  407-414.   2018.

\bibitem{xav14}G.B.A. Xavier, E. Suresh, and I. Gutman,    \textit{Counting relations for general Zagreb indices},      {Kragujevac J. Math.},     \textbf{38}, (1),   95-103. 2014. 

\bibitem{gut01}I. Gutman,  M. Lepovi\'{c},     \textit{Choosing the exponent in the definition of the connectivity index},     {J. Serb. Chem. Soc.},     \textbf{66}, (9),  605-611. 2001.

\bibitem{zho09}B. Zhou,  W. Luo,    \textit{A note on general Randi\'{c} index},   {MATCH Commun. Math. Comput. Chem.},    \textbf{62},   155-162.  2009.

\bibitem{cui17}Q. Cui,   L. Zhong,     \textit{The general Randi\'{c} index of trees with given number of pendent vertices},      \textbf{302}, (1),  111-121.  2017.  

\bibitem{alex14}V. Alexander,       \textit{Upper and lower bounds of symmetric division deg index},   {Iran. J. Math. Chem.},        \textbf{52},   91-98.  2014.

\bibitem{gup17}C.K. Gupta,  V. Lokesha,    B.S. Shwetha, and  P.S. Ranjini,     \textit{Graph operations on the symmetric division deg index of graphs},     {Palestine J. Math.},       \textbf{6}, (1),  280-286.  2017.

\bibitem{azar11}M. Azari,  A. Iranmanesh,      \textit{Generalized Zagreb index of graphs},    \textit{Studia Univ. Babes-Bolyai.},     \textbf{56}, (3), 59-70.  2011.

\bibitem{sar18}P. Sarkar,  N. De, and A. Pal,    \textit{The generalized Zagreb index of some carbon structures},         {Acta Chem. Iasi},   \textbf{26}, (1),  91-104.   2018.

\bibitem{sar19}P. Sarkar,  N. De,  I.N. Congul, and A. Pal,    \textit{The $(a,b)$-Zagreb index of some derived networks},         {J. Taibah Univ. Sci.},   \textbf{13}, (1),  79-86.  2019.

\bibitem{sark18}P. Sarkar,  N. De, and A. Pal,    \textit{The $(a,b)$-Zagreb index of nanostar dendrimers},         {U.P.B. Sci. Bull.},   \textbf{80}, (4),  67-82.  2018.

\bibitem{liu15}J.B. Liu,   X.F. Pan,      \textit{asymptotic incidence energy of lattices},   {Physica A},        \textbf{422},    193-202.   2015. 

\bibitem{chen16}L. Chen,  T. Li,   J. Liu,   Y. Shi, and  H. Wanhg,     \textit{On the Wiener polarity index of lattice networks},   {PLoS ONE},        \textbf{11}, (12),  1-28. 2016.



\end{thebibliography}
\end{document}